\documentclass[12pt] {article}
\usepackage{amsmath,amsthm}
\usepackage{amssymb,latexsym}
\title{Hard Lefschetz theorem for valuations, complex integral geometry, and
unitarily invariant valuations.}
\date{}
\author{Semyon Alesker
\\  { \normalsize Department of Mathematics, Tel Aviv University, Ramat Aviv}
 \\  { \normalsize 69978 Tel Aviv,
Israel }}
\newcommand{\RR}{\mbox{\rm $~\vrule height6.5pt width0.5pt
depth0.3pt\!\!$R}}

\newcommand{\CC}{\mbox{\rm $~\vrule height6.5pt width0.5pt
depth0.3pt\!\!$C}}

\newcommand{\HH}{\mbox{\rm $~\vrule height6.5pt width0.5pt
depth0.3pt\!\!$H}}
\newcommand{\DD}{\mbox{\rm $~\vrule height6.5pt width0.5pt
depth0.3pt\!\!$D}}

\def\eps{\varepsilon}

\def\lam{\lambda}

\def\str{\longrightarrow}

\def\ckl{C_{k,l}}
\def\gr{{}\!^ {\textbf{R}} Gr}
\def\grc{{}\!^ {\textbf{C}} Gr}
\def\agr{{}\!^ {\textbf{R}} {\cal A}Gr}
\def\agrc{{}\!^ {\textbf{C}} {\cal A}Gr}
\def\algr{ {\cal A}LGr}
\def\gro{\gr^+_{2,4}}
\def\val{Val^{ev}_{k}(V)}
\def\dgrc{\delta _{{}\!^{\textbf{C}} Gr}}
\def\tphi{\tilde \phi}
\def\tpsi{\tilde \psi}
\def\tf{\tilde f}
\def\tg{\tilde g}
\def\qed { q.e.d. }

\swapnumbers
\newtheorem{theorem}{Theorem}[subsection]
\newtheorem{corollary}[theorem]{Corollary}
\newtheorem{lemma}[theorem]{Lemma}
\newtheorem{proposition}[theorem]{Proposition}

\theoremstyle{definition}

\begin{document}
\maketitle
\begin{abstract}
We obtain new general results on the structure of the space of
translation invariant continuous valuations on convex sets (a
version of the hard Lefschetz theorem). Using these and our
previous results we obtain explicit characterization of unitarily
invariant translation invariant continuous valuations. It implies
new integral geometric formulas for real submanifolds in Hermitian
spaces generalizing the classical kinematic formulas in Euclidean
spaces due to Poincar\'e, Chern, Santal\'o, and others.
\end{abstract}
\setcounter{section}{-1}
\section{Introduction.}
\setcounter{subsection}{1} In this paper we obtain new results on
the structure of the space of even translation invariant
continuous valuations on convex sets. In particular we prove a
version of hard Lefschetz theorem for them and introduce certain
natural duality operator which establishes an isomorphism between
the space of such valuations on a linear space $V$ and on its dual
$V^*$ (with an appropriate twisting). Then we obtain an explicit
geometric classification of unitarily invariant translation
invariant continuous valuations on a Hermitian space $\CC^n$. This
classification is used to deduce new integral geometric formulas
for real submanifolds in Hermitian spaces generalizing the
classical kinematic formulas in Euclidean spaces due to
Poincar\'e, Chern, Santal\'o, and others.

Let us describe the results in more details. First let us remind
the definition of valuation. Let $V$ be a finite dimensional real
vector space. Let ${\cal K}(V)$ denote the class of all convex
compact subsets of $V$.

{\bf Definition.} {\itshape a) A function $\phi :{\cal K}(V) \str
\CC$ is called a valuation if for any $K_1, \, K_2 \in {\cal
K}(V)$ such that their union is also convex one has
$$\phi(K_1 \cup K_2)= \phi(K_1) +\phi(K_2) -\phi(K_1 \cap K_2).$$

b) A valuation $\phi$ is called continuous if it is continuous
with respect the Hausdorff metric on ${\cal K}(V)$. }

 Remind that
the Hausdorff metric $d_H$ on ${\cal K}(V)$ depends on the choice
of a Euclidean metric on $V$ and it is defined as follows:
$d_H(A,B):=\inf\{ \eps >0|A\subset (B)_\eps \mbox{ and } B\subset
(A)_\eps\},$ where $(U)_\eps$ denotes the $\eps$-neighborhood of a
set $U$. Then ${\cal K} (V)$ becomes a locally compact space (by
the Blaschke selection theorem).

In this paper we are interested only in translation invariant
continuous valuations. The space of such valuations will be
denoted by $Val(V)$.  The simplest examples of such valuations are
a Lebesgue measure on $V$ and the Euler characteristic $\chi$
(which is equal to 1 on each convex compact set). For the
classical theory of valuations we refer to the surveys
 \cite{mcmullen-survey}, \cite{mcmullen-schneider}. For a brief
overview of more recent results see \cite{alesker-ecm} and
\cite{alesker-icm}.

{\bf Definition.} {\itshape  A valuation $\phi$ is called
homogeneous of degree $k$ (or $k$-homogeneous) if for every convex
compact set $K$ and for every scalar $\lam >0$
$$\phi(\lam K)=\lam ^k \phi (K).$$
}

 Let us denote by $Val _k(V)$ the space of translation invariant
continuous $k$-homogeneous valuations.

{\bf Theorem.}(McMullen \cite{mcmullen-euler}){\itshape
$$Val(V)=\bigoplus_{k=0}^{n} Val_k(V),$$
where $n=\dim V$. }

 In particular note that the degree of
homogeneity is an integer between 0 and $n=\dim V$. It is known
that $Val_0(V)$ is one-dimensional and it is spanned by the Euler
characteristic $\chi$, and $Val_n(V)$ is also one-dimensional and
is spanned by a Lebesgue measure \cite{hadwiger-book}. The space
$Val_n(V)$ is also denoted by $|\wedge V^*| $ (or by $Dens(V)$,
the space of densities on $V$). Let us denote by $Val^{ev}(V)$ the
subspace of $Val(V)$ of even valuations (a valuation $\phi$ is
called even if $\phi(-K)=\phi(K)$ for every $K\in {\cal K}(V)$).
Similarly one defines the subspace $Val^{odd}(V)$ of odd
valuations. One has further decomposition with respect to parity:
$$Val_k(V)=Val_k^{ev}(V)\oplus Val_k^{odd}(V),$$
where $Val_k^{ev}(V)$ is the subspace of even $k$-homogeneous
valuations, and $Val_k^{odd}(V)$ is the subspace of odd
$k$-homogeneous valuations.

Let us fix on $V$ a Euclidean metric, and let $D$ denote the unit
Euclidean ball with respect to this metric.  Let us define on the
space of translation invariant continuous valuations an  operation
$\Lambda$ of mixing with the Euclidean ball $D$, namely
$$ (\Lambda \phi) (K):= \frac{d}{d\eps} \big|_{\eps =0} \phi(K +\eps D)$$
for any convex compact set $K$. Note that $\phi(K +\eps D)$ is a
polynomial in $\eps \geq 0$ by McMullen's theorem
\cite{mcmullen-euler}. It is easy to see that the operator
$\Lambda$ preserves parity and decreases the degree of homogeneity
by one. In particular we have
$$\Lambda: \val \str Val^{ev}_{ k-1}(V).$$
To formulate our first main result we will need one more
definition from the representation theory. Let $G$ be a Lie group.
Let $\rho$ be a continuous representation of $G$ in a Fr\'echet
space $F$. A vector $v\in F$ is called $G$-{\itshape smooth} if
the map $G\str F$ defined by $g\longmapsto g(v)$ is infinitely
differentiable. It is well known (and easy to prove) that smooth
vectors form a linear $G$-invariant subspace which is dense in
$F$. We will denote it by $F^{sm}$. It is well known (see e.g.
\cite{wallach}) that $F^{sm}$ has a natural structure of a
Fr\'echet space, and the representation of $G$ in $F^{sm}$ is
continuous with respect to this topology. In our situation the
Fr\'echet space $F=Val(V)$ with the topology of uniform
convergence on compact subsets of ${\cal K}(V)$, and $G=GL(V)$.
The action of $GL(V)$ on $Val(V)$ is the natural one, namely for
any $g\in GL(V),\, \phi\in Val(V)$ one has
$(g(\phi))(K)=\phi(g^{-1}K)$.

The following result is a version of the hard Lefschetz theorem.

{\bf Theorem 1.1.1.} {\itshape Let $n/2<k\leq n$. Then
$$\Lambda
^{2k-n}:(Val^{ev}_k(V))^{sm}\str(Val^{ev}_{n-k}(V))^{sm}$$ is an
isomorphism. In particular for $1\leq i\leq 2k-n$ the map
$$\Lambda
^{i}:(Val^{ev}_k(V))^{sm}\str(Val^{ev}_{k-i}(V))^{sm}$$ is
injective. }

 Our terminology is motivated by the classical hard Lefschetz
theorem (see e.g. \cite{griffiths-harris}) about the cohomology of
K\"ahler manifolds. To continue this analogy note that recently we
have observed \cite{alesker-mult} the natural multiplicative
structure on $(Val(V))^{sm}$ (see also \cite{alesker-icm}). More
precisely this space has natural structure of commutative
associative graded algebra (where the grading is given by the
degree of homogeneity). It satisfies a version on the Poincar\'e
duality with respect to these multiplication and grading.

The operator $\Lambda$ turns out to be closely related to so
called cosine transform on real Grassmannians, and the proof of
Theorem 1.1.1 is based on the solution of the cosine transform
problem by J. Bernstein and the author \cite{alesker-bernstein}
(some particular cases of this problem were solved previously by
Matheron \cite{matheron1} and Goodey, Howard, and Reeder
\cite{goodey-howard-reeder}).

Our next main result establishes connection between even
translation invariant continuous valuations on $V$ and on its dual
space $V^*$. In order to formulate it let us make an elementary
remark from linear algebra. Let $E\subset V$ be any
$k$-dimensional subspace. One has the canonical isomorphism
$|\wedge^n V|=|\wedge ^k E|\otimes |\wedge ^{n-k}(V/E)|$. Note
also that $V/E=(E^\perp)^*$. Hence we get the canonical
isomorphism
$$|\wedge ^k E^*|=|\wedge^{n-k}(E^\perp)^*|\otimes |\wedge^n
V^*|.$$ Then we have

{\bf Theorem 1.2.1.}{\itshape For any $k= 0,1,\dots ,n (=\dim V)$
there exists a natural isomorphism
$$\DD:\left( Val^{ev}_{k}(V)\right)^{sm} \tilde{\str}
\left( Val ^{ev}_{n-k}(V^*)\right)^{sm}\otimes |\wedge^n V^*|.$$
This isomorphism $\DD$ is uniquely characterized by the following
property: let $\phi \in Val^{ev}_k(V)$ and let $E\in Gr_k(V)$;
then $\phi|_E= \DD(\phi)|_{E^\perp}$ under the above
identification $|\wedge ^k E^*|=|\wedge^{n-k}(E^\perp)^*|\otimes
|\wedge^n V^*|.$}

 The proof of this theorem  uses the representation theoretical interpretation of the space
$Val^{ev}(V)$ given in \cite{alesker2}, where this space was
characterized as the unique irreducible submodule of some standard
$GL(V)$-module with smallest Gelfand-Kirillov dimension (of the
corresponding Harish-Chandra module).

Now let us discuss the translation invariant continuous valuations
invariant under some group $G$ of linear transformations of $V$.
This space will be denoted by $Val^G(V)$. If $G$ is the group of
orthogonal transformations $O(n)$ or special orthogonal
transformations $SO(n)$ the corresponding space of valuations is
described explicitly by the following famous result of H.
Hadwiger.

{\bf Theorem.} (Hadwiger, \cite{hadwiger-book}) {\itshape Let $V$
be $n$-dimensional Euclidean space. The intrinsic volumes $V_0,
V_1, \dots ,V_n$ form a basis of
$Val^{SO(n)}(V)(=Val^{O(n)}(V))$.}

Let us remind the definition of the intrinsic volumes $V_i$. Let
$\Omega$ be a compact (not necessarily convex) domain in a
Euclidean space $V$ with  smooth boundary $\partial \Omega$. Let
$n=\dim V$. For any point $s\in \partial \Omega$ let $k_1(s),
\dots, k_{n-1}(s)$ denote the principal curvatures at $s$. For
$0\leq i \leq n-1$ define
$$V_i(\Omega):=\frac{1}{(n-i)vol_{n-i}(D_{n-i})} {n-1\choose n-1-i}^{-1}
\int _{\partial \Omega}\{k_{j_1}, \dots,k_{j_{n-1-i}}\} d\sigma,
$$ where $\{k_{j_1}, \dots,k_{j_{n-1-i}}\}$ denotes the $(n-1-i)$-th
elementary symmetric polynomial in the principal curvatures,
$d\sigma$ is the measure induced on $\partial \Omega$ by the
Euclidean metric, and $D_{n-i}$ denotes the unit
$(n-i)$-dimensional ball. It is well known (see
e.g.\cite{schneider-book}) that $V_i$ (uniquely) extends by
continuity in the Hausdorff metric to ${\cal K}(V)$. Define also
$V_n(\Omega):=vol(\Omega)$. Note that $V_0$ is proportional to the
Euler characteristic $\chi$. It is well known that $V_0, V_1,
\dots ,V_n$  belong to $Val^{O(n)}(V)$. It is easy to see that
$V_k$ is homogeneous of degree $k$.

Now let us describe unitarily invariant valuations on the
Hermitian space $\CC^n$. Let us denote by $IU(n)$ the group of
isometries of the Hermitian space $\CC^n$ preserving the complex
structure (then $IU(n)=\CC^n \rtimes U(n)$). Let $\agrc _{j} $
denote the Grassmannian of affine complex subspaces of $\CC^n$ of
complex dimension $j$. Clearly $\agrc _{j} $ is a homogeneous
space of $IU(n)$ and it has a unique (up to a constant)
$IU(n)$-invariant measure (called Haar measure). For every
non-negative integers $p$ and $k$ such that $2p\leq k\leq 2n$  let
us introduce the following valuations:
$$U_{k,p}(K)=\int_{E\in \agrc _{n-p}} V_{k-2p}(K\cap E) \cdot
dE.$$ Then $U_{k,p}\in Val_k^{U(n)}(\CC^n)$.

{\bf Theorem 2.1.1.}{\itshape The valuations $U_{k,p}$ with $0\leq
p \leq \frac{min\{k, 2n-k\}}{2}$ form a basis of the space
$Val_k^{U(n)}(\CC^n)$. }

This result is the Hermitian generalization of the (Euclidean)
Hadwiger theorem. The proof of this theorem is highly indirect. It
turns out to be necessary to study the $GL_{2n}(\RR)$-module
structure of the infinite dimensional space $Val^{ev}(\CC^n)$. The
proof of Theorem 2.1.1 uses most of the facts known about even
valuations including the solution of the McMullen conjecture
\cite{alesker2}, cosine transform \cite{alesker-bernstein}, the
hard Lefschetz theorem for valuations, and the results of Howe and
Lee \cite{howe-lee} on the $K$-type structure of certain
$GL$-modules.

Note that there are some other natural examples of valuations from
$Val^{U(n)}(\CC^n)$, for instance the averaged volume of
projections of a convex set to all complex (or, say, Lagrangian)
subspaces. Theorem 2.1.1 implies that all of them are linear
combinations of $U_{k,p}$ with the above range of indices $k,p$.
We would also like to mention another interesting example of such
valuation which comes from the complex analysis. It is so called
Kazarnovskii's pseudovolume. It was introduced and studied by B.
Kazarnovskii \cite{kazarnovskii1}, \cite{kazarnovskii2} in order
to write down a formula for the number of zeros  of a system of
exponential sums in terms of their Newton polytopes. His results
generalize in some sense the well known results of D. Bernstein
\cite{bernstein} and A. Kouchnirenko \cite{kouchnirenko} on the
number of zeros of a system of polynomial equations (see also
\cite{gelfand-kapranov-zelevinsky}). We will recall the definition
of Kazarnovskii's pseudovolume in Subsection 3.3. As a corollary
of Theorem 2.1.1 we present a new formula for Kazarnovskii's
pseudovolume in integral geometric terms (Theorem 3.3.2). It also
seems that the valuation property of Kazarnovskii's pseudovolume
was not mentioned previously in the literature.

The classification of unitarily invariant valuations is used to
obtain new integral geometric formulas in the Hermitian space
$\CC^n$. Let us state some of them.  Let $\Omega_1,\, \Omega_2$ be
compact domains with smooth boundary in $\CC^n$ such that
$\Omega_1 \cap U(\Omega_2)$ has finitely many components for all
$U\in IU(n)$. The new result is

{\bf Theorem 3.1.1.} {\itshape Let $\Omega_1, \, \Omega_2$ be
compact domains in $\CC^n$ with piecewise smooth boundaries such
that for every $U\in IU(n)$ the intersection $\Omega_1\cap
U(\Omega _2)$ has finitely many components. Then
$$\int_{U\in IU(n)} \chi (\Omega_1\cap U(\Omega_2)) dU=
\sum_{k_1+k_2=2n}\sum_{p_1, p_2} \kappa(k_1,k_2,p_1,p_2)
U_{k_1,p_1}(\Omega_1) U_{k_2,p_2}(\Omega_2),$$ where the inner sum
runs over $0\leq p_i\leq k_i/2, \, i=1,2$, and
$\kappa(k_1,k_2,p_1,p_2)$ are certain uniquely defined constants
depending on $n,k_1,k_2,p_1,p_2$ only. }

The study of the left hand side in this formulas was started by J.
Fu \cite{fu}.

{\bf Theorem 3.1.2.}{\itshape Let $\Omega$ be a compact domain in
$\CC^n$ with piecewise smooth boundary. Let $0<q<n, \, 0<2p<k<2q$.
Then
$$\int_{E\in \agrc_q} U_{k,p}(\Omega \cap E)=\sum
_{p=0}^{[k/2]+n-q}\gamma_p \cdot U_{k+2(n-q),p}(\Omega),$$ where
the constants $\gamma_p$ depend only on $n,\,q$, and $p$.}

Let us denote by $\algr(\CC^n)$ the (non-compact) Grassmannian of
{\itshape affine} Lagrangian subspaces of $\CC^n$. Clearly it is a
homogeneous space of the group $IU(n)$ and hence has a Haar
measure.

{\bf Theorem 3.1.3.} {\itshape Let $\Omega$ be a compact domain in
$\CC^n$ with piecewise smooth boundary. Then
$$\int_{\algr (\CC^n)} \chi(E\cap\Omega)dE=
\sum_{p=0}^{[n/2]}\beta_p\cdot U_{n,p}(\Omega),$$ where $\beta_p$
are certain uniquely defined constants depending on $n$ and $p$
only. }

Theorems 3.1.1 and 3.1.2 are analogs of general kinematic formulas
of Chern \cite{chern1}, \cite{chern2} and Federer \cite{federer}
(see also \cite{santalo}, especially Ch. 15, and
\cite{klain-rota}). Further generalizations in the Euclidean case
were obtained by Cheeger, M\"uller, and Schrader
\cite{cheeger-muller-schrader} and J. Fu \cite{fu}. For more
recent results in this direction and further references we refer
to the recent survey by Hug and Schneider \cite{hug-schneider}.
For classical results in Hermitian integral geometry  we refer to
\cite{chern3}, \cite{griffiths}, \cite{santalo-paper}. In these
papers the authors discuss the integral geometry of {\itshape
complex} submanifolds. The integral geometry of Lagrangian
submanifolds also was studied (see e.g. \cite{le-hong-van}). The
integral geometry of real submanifolds in the complex projective
space $\CC P^n$ was studies by H. Tasaki \cite{tasaki1},
\cite{tasaki2} and  Kang and Tasaki \cite{kang-tasaki1},
\cite{kang-tasaki2}. In these papers the authors obtain explicit
Poincar\'e type formulas for real submanifolds of certain specific
dimensions. Their results use in turn a general Poincar\'e type
formula in Riemannian homogeneous spaces due to R. Howard
\cite{howard}.

It would be of interest to compute the constants
$\kappa(k_1,k_2,p_1,p_2)$, $\gamma_p$, and $\beta_p$ in Theorems
3.1.1, 3.1.2, and 3.1.3 explicitly. We could not do it in general.
But we were able to compute them only in the first non-trivial
case $n=2$. One of such computations is presented in Subsection
3.2.  Much more complete treatment of the integral geometric
formulas (including computation of all constants) in $\CC^2$,
$\CC^3$, and also in the 2- and 3- dimensional complex projective
and hyperbolic spaces was done recently by H. Park in his thesis
\cite{park}.

The paper is organized as follows. In Section 1 we discuss the
results about the structure of the space of even translation
invariant continuous valuations. Namely in Subsection 1.1 we prove
the hard Lefschetz theorem for valuations and deduce some
corollaries from it. In Subsection 1.2 we discuss the duality on
valuations, in particular we prove Theorem 1.2.1. In Section 2 we
prove the classification of unitarily invariant translation
invariant continuous valuations. In Section 3 we discuss the
integral geometry in complex spaces. In Subsection 3.1 we obtain
the integral geometric formulas in $\CC^n$. In Subsection 3.2 we
compute explicitly the constants in one of such formulas in
$\CC^2$ (Theorem 3.2.4).

{\bf Acknowledgements.} We express our gratitude to J. Bernstein
for numerous useful discussions. We would also like to thank P.
Biran, L. Polterovich, and J. Fu for useful conversations.

\section{Hard Lefschetz theorem and duality for valuations.}
Let $V$ be an $n$-dimensional real vector space. In  Subsection
1.1 of this section we prove an analogue of the hard Lefschetz
theorem for translation invariant even continuous valuations. In
Subsection 1.2 we introduce the notion of a valuation dual to a
given translation invariant even continuous valuation which
satisfies some additional mild technical condition of $GL(V)$-
smoothness (defined in the introduction). This construction uses
the representation theoretical interpretation of the space of
valuations given in \cite{alesker2}. The geometric examples will
be given in Proposition 2.1.7 of Section 2.
\subsection{An analogue of the hard Lefschetz theorem
for valuations.} The main result of this subsection is the
following analogue of the hard Lefschetz theorem where the
operator $\Lambda$ was defined in the introduction.
\begin{theorem}
Let $n\geq k>n/2$. Then $\Lambda ^{2k-n} : (\val)^{sm} \str
(Val^{ev}_{n-k}(V))^{sm}$ is an isomorphism. In particular
$\Lambda ^i : \val \str Val^{ev}_{k-i}(V)$ is injective for $1\leq
i \leq 2n-k$.
\end{theorem}

The proof of this theorem uses the cosine transform on real
Grassmannians, thus we will remind first its definition and the
relevant properties. We will denote by R$\gr_j (V)$ the
Grassmannian of real $j$-dimensional subspaces in $V$. Assume that
$1\leq i\leq j \leq n-1$. For two subspaces $E\in \gr_i(V)$, $F\in
\gr_j(V)$ let us define the {\itshape cosine of the angle} between
$E$ and $F$:
$$|cos(E,F)|:=\frac{vol_i(Pr_F(A))}{vol_i(A)},$$
where $A$ is any subset of $E$ of non-zero volume, $Pr_F$ denotes
the orthogonal projection onto $F$, and $vol_i$ is the
$i$-dimensional measure induced by the Euclidean metric. (Note
that this definition does not depend on the choice of a subset
$A\subset E$). In the case $i\geq j$ we define the cosine of the
angle between $E$ and $F$ as cosine of the angle between their
orthogonal complements:
$$|cos(E,F)| :=|cos(E^\perp ,F^\perp )|.$$
(It is easy to see that if $i=j$ both definitions coincide.)

For any $1\leq i, \, j \leq n-1$ one defines the cosine transform
$$T_{j,i}:C(\gr_i(V)) \str C(\gr_j(V))$$ as follows:
$$(T_{j,i}f)(F):= \int_{\gr_i(V)} |cos(E,F)| f(E) dE,$$
where the integration is with respect to the Haar measure on the
Grassmannian such that the total measure is equal to 1. Clearly
the cosine transform commutes with the action of the orthogonal
group $O(n)$, and hence its image is an $O(n)$-invariant subspace
of functions.

Now let us recall the imbedding $\val \str C(\gr_k(V))$ which we
will call the Klain imbedding. Let $\phi \in \val$. For every  $E
\in \gr_k(V)$ let us consider the restriction of $\phi$ to all
convex compact subsets of $E$. This is an even translation
invariant valuation homogeneous of degree $k$. Hence, by a result
due to Hadwiger \cite{hadwiger-book}, it is a density on $E$ (i.e.
a Lebesgue measure). Thus it is equal to $f(E) \cdot vol_E$, where
$vol_E$ is the volume form on $E$ defined by the metric on $V$,
and $f(E)$ is a constant depending on $E$. Thus $\phi \mapsto f$
defines the map $\val \str C(\gr_k(V))$ which turns out to be an
imbedding by a result due to D. Klain (\cite{klain1}; this result
was stated in this form in \cite{klain2} and in \cite{alesker1}).
Let us denote this image by $I_k$. Moreover it was shown in
\cite{alesker-bernstein} that the image of the Klain imbedding
coincides with the image of the cosine transform
$T_{k,k}:C(\gr_k(V)) \str C(\gr_k(V))$ (at least on the level of
$GL(V)$-smooth vectors).
\begin{lemma}
Let $k\geq n/2$. The cosine transform
$$T_{n-k,k}:C(\gr_k(V)) \str C(\gr_{n-k}(V))$$
maps $I_k$ to $I_{n-k}$ and induces isomorphism of $O(n)$-smooth
vectors of these subspaces.
\end{lemma}
{\bf Proof.} It is well known that for admissible $GL(V)$-modules
of finite length the subspaces of $GL(V)$-smooth and $O(n)$-smooth
vectors coincide (more generally, $GL(V)$ can be replaced by any
real reductive group $G$, and $O(n)$ can be replaced by a maximal
compact subgroup of $G$). First let us prove that $I_k$ and
$I_{n-k}$ have the same decomposition under the action of the
orthogonal group $O(n)$. Indeed the correspondence $E\mapsto
E^\perp$ induces an isomorphism $S:C(\gr_k(V)) \str
C(\gr_{n-k}(V))$ commuting with the action of $O(n)$. Moreover we
have the following relation between the cosine transforms:
$$T_{n-k,n-k}=S T_{k,k} S^{-1}.$$
Hence it follows that $S((I_k)^{sm})=(I_{n-k})^{sm}$ (it is
immediate on the level of $O(n)$-finite vectors; to deduce it for
$O(n)$-smooth vectors one should use the Casselman-Wallach theorem
\cite{casselman} as it is done in \cite{alesker-bernstein}).

Next it is well known (see e.g. \cite{alesker-bernstein}, Lemma
1.7) that the cosine transform $T_{n-k,k}$ can be written (up to a
non-zero normalizing constant which we ignore) as a composition
$T_{n-k,n-k}\circ R_{n-k,k}$, where $R_{n-k,k}:C(\gr_k(V)) \str
C(\gr_{n-k}(V))$ is the Radon transform. It was shown in
\cite{gelfand-graev-rosu} that
$$R_{n-k,k}:C^\infty(\gr_k(V)) \str C^\infty(\gr_{n-k}(V))$$ is an isomorphism. We claim that
$R_{n-k,k}((I_k)^{sm})=(I_{n-k})^{sm}$. To see this remind that
the quasiregular representation of $O(n)$ in the space of
functions on the Grassmannians is multiplicity free (since the
Grassmannians are symmetric spaces). Hence it follows that two
$O(n)$-invariant subspaces of $C(\gr_{n-k}(V))$ have the same
$O(n)$-finite vectors if and only if these subspaces have the same
decomposition under the action of $O(n)$ (in the abstract sense).
Hence $R_{n-k,k}(I_k)$ and $I_{n-k}$ have the same $O(n)$-finite
vectors. The coincidence of $O(n)$-smooth vectors follows again
from the Casselman-Wallach theorem \cite{casselman} and the fact
that the Radon transform can be rewritten as an intertwining
operator of admissible $GL(V)$-modules of finite length (see
\cite{gelfand-graev-rosu}).

Since $T_{n-k,n-k}$ is selfadjoint its restriction to $I_{n-k}$
has trivial kernel and dense image. But the key observation of
\cite{alesker-bernstein} was that $T_{n-k,n-k}$ can be rewritten
as an intertwining operator of certain $GL(V)$-modules. This and
the Casselman-Wallach theorem \cite{casselman} imply that
$T_{n-k,n-k}((I_{n-k})^{sm})=(I_{n-k})^{sm}$. Hence
$T_{n-k,k}((I_{k})^{sm})=(I_{n-k})^{sm}$. \qed

 Now let us prove Theorem 1.1.1.

{\bf Proof} of Theorem 1.1.1. Since the image of $GL(V)$-smooth
continuous $k$-homogeneous valuations in $C(\gr_k(V))$ coincides
with the image of the cosine transform on $GL(V)$- smooth
functions, then every $GL(V)$-smooth valuation $\phi \in \val$ can
be represented in the form
$$\phi(K)= \int_{\gr_k(V)} f(E) vol _k (Pr_E(K)) dE,$$
where $f$ is a smooth function on $\gr_k(V)$, $K$ is an arbitrary
convex compact set, $Pr_E$ denotes the orthogonal projection onto
$E$, and the integration is with respect to the Haar measure on
the Grassmannian. Moreover for every smooth function $f$, the
expression defined by this formula is a valuation from
$(\val)^{sm}$. For a given valuation $\phi$ the function $f$ is
not defined uniquely. But we can choose $f\in I_k$, i.e. in the
image of the cosine transform;
 then it will be defined uniquely. So we will
assume that $f\in I_k$.
Let us apply $\Lambda ^{2k-n}$ to it. Then it is easy to see that
$$(\Lambda ^{2k-n} \phi)(K)= c \cdot \int_{\gr_k(V)} f(E) V_{n-k} (Pr_E(K)) dE, $$
where $c$ is a non-zero normalizing constant, and
$V_{n-k}(Pr_E(K)) $ denotes the $(n-k)$-th intrinsic volume of
$Pr_E(K)$ inside $E$, i.e. it is the mixed volume of $Pr_E(K)$
taken $n-k$ times with the unit ball of $E$ taken $2k-n$ times.
The image $g$ of $\Lambda ^{2k-n} \phi $ in functions on the
Grassmannian $C(\gr_{n-k}(V))$ can be described as follows. It is
easy to see that for every subspace $F\in \gr_{n-k}(V)$
$$g(F)= c' \cdot \int_{\gr_{n-k}(V)} f(E) |cos(F,E)| dE,$$
where $c'$ is a non-zero normalizing constant. Namely $g$ is equal
(up to a normalization) to the cosine transform $T_{n-k,k}(f)$ of
$f$. By Lemma 1.1.2 $T_{n-k,k}$ induces the isomorphism between
$GL(V)$-smooth vectors of $I_k$ and of $I_{n-k}$. This proves
Theorem 1.1.1.
 \qed

For a subgroup $G \subset GL(V)$ let us denote by $Val^G_k(V)$ the
space of translation invariant $G$- invariant $k$- homogeneous
continuous valuations. Let $h_k:=\dim Val_k^G(V)$.
\begin{corollary}
Let $G$ be a compact subgroup of the orthogonal group which acts
transitively on the unit sphere and contains the operator $-Id$.
Then $Val ^G_{k}(V)$ is a finite dimensional space, and  for $n/2
< k \leq n$
$$ \Lambda ^{2k-n}: Val ^G_{k}(V) \str Val^G_{n-k}(V)$$
is an isomorphism. Consequently the numbers $h_i$ satisfy
the Lefschetz inequalities:
$$h_i \leq h_{i+1} \mbox{ for } i <n/2, \mbox{ and } h_i=h_{n-i} \mbox{ for } i=0, \dots ,n.$$
\end{corollary}
{\bf Proof.} The finite dimensionality of $Val ^G_{k}(V)$ was
proved in \cite{alesker1}. Let us show that this implies that all
vectors from $Val ^G_{k}(V)$ are $O(n)$-finite (in particular
$GL(V)$- smooth). Indeed let $Z$ be the minimal closed $O(n)$-
invariant subspace of the space $Val_{k}(V)$  containing $Val
^G_{k}(V)$. The space $Z$ is decomposed under the action of $O(n)$
into the direct sum of irreducible components, and each component
enters with finite multiplicity (since the space of translation
invariant continuous valuations of the given degree of homogeneity
and parity can be realized as a subquotient of a representation of
$GL(V)$ induced from a character of a parabolic subgroup, see
Section 2 in \cite{alesker2}). Thus let $Z=\oplus _i \rho _i$ be
this decomposition. We have a continuous projection $\pi:
Val_{k}(V) \str Val^G_k(V)$ defined by $\pi (\phi)= \int_{g \in G}
g(\phi) dg$. Clearly $Im (\pi)=Val^G_k(V) =(Z)^G$. But
$(Z)^G=\oplus _i (\rho _i)^G$. Since $Val^G_k(V)$ is finite
dimensional, $(\rho _i)^G=0$ for all but finitely many $i$'s. In
other words there is a finite set of indices $A$ such that $Val
^G_{k}(V) \subset \oplus_{i\in A} \rho_i$. Thus all elements of
$Val ^G_{k}(V)$ are $O(n)$-finite.

Next obviously $\Lambda (Val^G_k(V)) \subset Val^G_{k-1}(V) $. The
rest follows from Theorem 1.1.1.
 \qed

\subsection{Duality on valuations.} Let $V$ be an $n$-dimensional
real vector space. Let us denote by $V^*$ its dual space. Let us
denote by $|\wedge^n V^*|$ the (one-dimensional) space of
complex-valued Lebesgue measures on $V$. Let us consider the space
$Val^{ev}_{k}(V^*)\otimes |\wedge^n V^*|$ of translation invariant
even continuous $k$-homogeneous valuations on $V^*$  with values
in $|\wedge^n V^*|$. Note that on both spaces we have the natural
(continuous) representation of the group $GL(V)$.

Before we state the main result of this subsection let us make a
remark. For any subspace $E\in Gr_k(V)$ consider the short exact
sequence $0\str E\str V\str V/E\str 0$. From this sequence one
gets the canonical isomorphism $|\wedge^n V|=|\wedge ^k E|\otimes
|\wedge ^{n-k}(V/E)|$. Note also that $V/E=(E^\perp)^*$. Hence we
get the canonical isomorphism
$$|\wedge ^k E^*|=|\wedge^{n-k}(E^\perp)^*|\otimes |\wedge^n
V^*|.$$
 The main result of this subsection is
\begin{theorem}
For any $k= 0,1,\dots ,n$ there exists a natural isomorphism
$$\DD:\left( Val^{ev}_{k}(V)\right)^{sm} \tilde{\str}
\left( Val ^{ev}_{n-k}(V^*)\right)^{sm}\otimes |\wedge^n V^*|.$$
This isomorphism $\DD$ is defined uniquely by the following
property: let $\phi \in Val^{ev}_k(V)$ and let $E\in Gr_k(V)$;
then $\phi|_E= \DD(\phi)|_{E^\perp}$ under the above
identification $|\wedge ^k E^*|=|\wedge^{n-k}(E^\perp)^*|\otimes
|\wedge^n V^*|.$
\end{theorem}
{\bf Proof.} First let us rewrite the Klain imbedding we have
discussed in Subsection 1.1 of the space of even valuations $Val
^{ev}_k(V)$ to functions on the Grassmannian in the notation which
does not use any Euclidean structure. Instead of functions on the
Grassmannian we have to consider sections of certain line bundle
$L_k$ over the Grassmannian $\gr_k(V)$. The fiber of $L_k$ over a
subspace $E\in \gr_k(V)$ is the (one-dimensional) space $|\wedge^k
E^*|$ of complex valued Lebesgue measures on $E$. Clearly $L_k$ is
naturally $GL(V)$- equivariant. Let us denote by $C(\gr_k, L_k)$
the space of continuous sections of $L_k$. The map we have
described in Subsection 1.1 can be rewritten as follows. Fix a
valuation $\phi \in Val^{ev}_{k}(V)$. For any $E\in \gr_k(V)$ let
us consider the restriction of $\phi$ to $E$. As previously, since
this restriction $\phi |_E$ has maximal degree of homogeneity
(equal to $k$) by Hadwiger's theorem \cite{hadwiger-book} $\phi
|_E$ is a Lebesgue measure on $E$. Thus $\phi$ defines a
continuous section of $L_k$. As we have mentioned, the constructed
map is injective. One of the main results of \cite{alesker2} says
that the image of $Val^{ev}_{k}(V)$ in $C(\gr_k (V), L_k)$ under
this map is the unique "small" irreducible $GL(V)$-submodule
(Theorem 1.3 combined with Theorem 3.1 in \cite{alesker2}) . Let
us give some comments what does it mean "small". First replace all
$GL(V)$- modules by their Harish-Chandra modules which are purely
algebraic objects. For each Harish-Chandra module one defines an
{\itshape associated variety} (or Bernstein's variety) which is an
algebraic subvariety of the Lie algebra $gl_n(\CC)$, where
$n=dim(V)$. For details we refer to \cite{borho-brylinski}. When
we say that a given irreducible submodule $A$ of a module $B$ is
"small" it means that the dimensions of the associated varieties
of all other irreducible subquotients of $B$ are strictly greater
than that of $A$. Note also that the dimension of the associated
variety of $A$ is equal to the Gelfand-Kirillov dimension of the
underlying Harish-Chandra module.

Now let us continue constructing the isomorphism $\DD$. Let us
consider the line bundle $M_k$ over $\gr _{n-k}(V^*)$ the fiber of
which over any $F\in \gr _{n-k}(V^*)$ is equal to $|\wedge^{n-k}
F^*|\otimes |\wedge ^n V^*|$ (note that $|\wedge^{n-k} F^*|$ is
identified with the space of Lebesgue measures on $F$). As
previously, $Val ^{ev}_{n-k}(V^*)\otimes |\wedge^n V^*|$ can be
realized as the only "small" irreducible submodule of $C(\gr
_{n-k}(V^*), M_k)$ (indeed these spaces differ from the previous
two only by the twist by $|\wedge^n V^*|$). Hence it is sufficient
to present the natural isomorphism between $C^{\infty}(\gr _k(V),
L_k)$ and $C^{\infty}(\gr _{n-k}(V^*), M_k)$ where $C^{\infty}$
denotes the space of $C^{\infty}$- sections of the bundles. Let us
do it. Let $E\in \gr_k(V)$. As previously,  we have the canonical
isomorphism
$$|\wedge ^k
E^*|=|\wedge ^{n-k} (E^\perp)^* |\otimes |\wedge ^n V^*|.$$ The
correspondence $E \longmapsto E^\perp$ and the last identification
give the desired isomorphism. \qed

Now let us assume that $V$ is a Euclidean space, i.e. on $V$ we
are given a positive definite quadratic form. This gives us the
identification of $V$ with its dual space $V^*$, and the
identification of the space $|\wedge^n V^*|$ of Lebesgue measures
on $V$ with the complex line $\CC$ (such that $1\in \CC$
corresponds to the Lebesgue measure on $V$ which is equal to 1 on
the unit cube). Also for any subspace $E$ let us denote by $vol_E$
the Lebesgue measure on $E$ which is equal to 1 on the unit cube.
 Under these identifications we get
$$\DD: (Val^{ev}_{k}(V))^{sm}\tilde{ \str} (Val^{ev}_{n-k}(V))^{sm}.$$
For this operator we have the following result which can be easily
deduced from the last theorem.
\begin{theorem}
Let $V$ be an $n$-dimensional Euclidean space. Then for any
$k=0,1,\dots ,n$
$$\DD: (Val^{ev}_{k}(V))^{sm}\tilde{ \str} (Val^{ev}_{n-k}(V))^{sm}$$
is an isomorphism and $\DD ^2=Id$. This operator $\DD$ is defined
uniquely by the following property: let $\phi \in Val^{ev}_k(V)$
and let $E\in Gr_k(V)$; if $\phi|_E =f(E)\cdot vol_E$ then $\DD
\phi |_{E^\perp} = f(E) vol_{E^\perp}$. Also $\DD$ commutes with
the action of $O(n)$.
\end{theorem}
{\bf Example.} Let $\chi$ denote the Euler characteristic on a
Euclidean space $V$. Clearly $\chi \in Val_0(V)$. Then $\DD(\chi)=
vol_V$, and $\DD(vol_V) =\chi$.

\section{ Unitarily invariant valuations.}
\setcounter{subsection}{1} \setcounter{theorem}{0}
 In this section we will describe
 unitarily invariant translation invariant continuous valuations
on convex compact subsets of $\CC^n$ by writing down explicitly a
basis in this space.
Let $k, \, l$ be integers such that $0\leq k \leq 2n$ and  $k/2 \leq l \leq n$.
Let us define a valuation
$$\ckl (K):=\int _{\grc _{l,n}} V_k (Pr_F (K)) dF,$$
where the integration is with respect to the Haar measure on the
complex Grassmannian of complex $l$-dimensional subspaces in $\CC
^n$, $Pr_F$ denotes the orthogonal projection onto $F$, and $V_k
(Pr_F (K))$ denotes the $k$-th intrinsic volume of $Pr_F (K)$
inside $F$, namely it is the mixed volume of $Pr_F (K)$ taken $k$
times with the unit Euclidean ball in $F$ taken $2l-k$ times.
Clearly $\ckl \in Val^{U(n)}_k(\CC^n)$. Note that for $l=n$ we get
the usual intrinsic volumes. For $k=0$ we get the Euler
characteristic, and for $k=2n, \, l=n$ we get the Lebesgue
measure. Our next main result is
\begin{theorem}
Let $k$ be an integer, $0\leq k \leq 2n$. The dimension of the
space $Val_{k}^{U(n)}(\CC^n)$ is equal to $1+ min\{[k/2], \,
[(2n-k)/2] \}$.
 The valuations $\ckl$ with $\frac{ max \{k, 2n-k\}}{2} \leq l \leq n$,
form a basis of $Val_{k}^{U(n)}(\CC^n)$.
\end{theorem}

{\bf Remark.} Later on in this section we will present another
basis in the space of unitarily invariant valuations. This basis
will be more convenient for the applications in integral geometry
and for non-convex sets. In fact the connection between these two
bases is not quite trivial and leads to new integral geometric
formulas. This material will be discussed in more detail in
Section 3.

{\bf Proof.} The dimension of $Val_{k}^{U(n)}(\CC^n)$ was computed
in \cite{alesker2}. Hence it remains to show that the valuations
$\ckl $ with $\frac{ max \{k, 2n-k\}}{2} \leq l \leq n$ are
linearly independent. First of all it is clear that $C_{k,l}= c
\cdot \Lambda (C_{k+1,l}) $, where $\Lambda$ is the operator from
the hard Lefschetz theorem (Theorem 1.1.1), and $c$ is a non-zero
constant depending on $n,k,l$ only. Hence by the hard Lefschetz
theorem for unitarily invariant valuations (Corollary 1.1.3) the
statement is reduced to the case $k \geq n$. Let us prove this
case. We will prove the statement by induction in $2n-k$. If
$2n-k=0$ then the result is clear since by \cite{hadwiger-book}
any translation invariant continuous $N$-homogeneous valuation on
$\RR ^N$ is a Lebesgue measure. Now assume that $n\leq k <2n$, and
the theorem is true for valuations homogeneous of degree $>k$. If
$k$ is odd then the induction assumption, Corollary 1.1.3, and the
computation of the dimension of unitarily invariant
$k$-homogeneous valuations imply the result. Hence let us assume
that $k$ is even. Again using Corollary 1.1.3 it is sufficient to
check that $C_{k,k/2}$ can not be presented as a linear
combination of valuations $C_{k,l}$ with $l> \frac{k}{2}$.

In order to prove it, we will show that the special orthogonal
group $SO(2n)$ acts differently on $C_{k,k/2}$ and on $ C_{k,l}$
with $l> \frac{k}{2}$. To formulate this more precisely let us
introduce some notation. First recall that the set of highest
weights of $SO(2n)$ is parameterized by sequences of integers
$\mu_1, \dots , \mu_{n-1}, \mu_n$ such that $\mu_1 \geq \dots \geq
\mu_{n-1}\geq |\mu_n|$.
 For $1\leq l \leq n$
let us denote by $\Lambda (l)$ the subset of highest weights of
$SO(2n)$ such that all $\mu_i$'s are even and if $l<n$ satisfy in
addition the following condition: $\mu_j=0 $ for $j>l$.

The following result was proved in \cite{alesker2}, Proposition
6.3.
\begin{lemma}
The natural representation of $SO(2n)$ in the space
$Val^{ev}_{k}(\CC^n)$ is multiplicity free and is isomorphic to a
direct sum of irreducible components with highest weights $(\mu_1,
\mu_2, \dots , \mu_n) \in \Lambda (min(k,2n-k))$ such that
$|\mu_2| \leq 2$.
\end{lemma}
Note that the explicit description of the $K$-type structure was
heavily based on the results of Howe and Lee \cite{howe-lee}.

The next result is well known (see e.g. \cite{takeuchi}, \S 8).
\begin{lemma}
In every irreducible representation of $SO(2n)$ the subspace of
$U(n)$-invariant vectors is at most 1-dimensional. This subspace
is 1-dimensional if and only if the highest weight of the
irreducible representation of $SO(2n)$ is of the form $(\mu_1,
\dots ,\mu_n)$ where

(i) if $n$ is even then
$$\mu_1 =\mu_2 \geq \mu_3 =\mu_4 \geq \dots \geq \mu_{n-1}=\mu_n
\geq 0;$$

(ii) if $n$ is odd then
$$\mu_1 =\mu_2 \geq \mu_3 =\mu_4 \geq \dots \geq \mu_{n-2}=
\mu_{n-1} \geq \mu_n =0.$$
\end{lemma}

 The following lemma and Corollary 1.1.3 obviously imply Theorem 2.1.1.
\begin{lemma} Let $k$ be even, $n \leq k <2n$.
(i) The valuations $C_{k,l}$ with $l>k/2$ belong to the sum of the
representations with highest weights $\mu \in \Lambda (2n-k-1)$.

(ii) The valuation $C_{k,k/2}$ does not belong to the above sum.
\end{lemma}

{\bf Proof.} First let us prove part (i) of the lemma. As we have
mentioned earlier $C_{k,l} =c \cdot \Lambda (C_{k+1,l})$ if
$l>k/2$. Since the operator $\Lambda$ commutes with the action of
$SO(2n)$ on valuations then it is sufficient to check that
$C_{k+1,l}$ belongs to the sum of irreducible components with
highest weights from $\Lambda (2n-k-1)$. As it was mentioned in
Section 1 of this paper the space $Val^{ev}_{k+1}(\CC^n)$ can be
imbedded into the space of continuous functions $C(\gr_{k+1,2n})$.
But it is well known that all irreducible representations of
$SO(2n)$ which appear in the last space belong to $\Lambda
(2n-k-1)$ (see e.g. \cite{takeuchi} \S 8 for the general case of
compact symmetric spaces). This proves part (i) of the lemma.

Let us prove part (ii) which is somewhat more computational. We
will show that the image of $C_{k,k/2}$ in $C(\gr_{k,2n})$ is not
orthogonal to the irreducible subspace in $C(\gr_{k,2n})$ with
highest weight $(\underbrace{2,2, \dots, 2}_{2n-k \mbox{ times }},
0,\dots , 0)$. Clearly this will finish the proof of Lemma 2.1.4,
and hence the proof of Theorem 2.1.1.

From the definition of $C_{k,k/2}$ we immediately see that its
image in $C(\gr_{k,2n})$ is the function $f$ such that
$$ f(E)= c \cdot \int _{\grc _{k/2,n}} |cos(E,F)| dF ,$$
where $c$ is a non-zero normalizing constant. In other words $f$
is proportional to the cosine transform of the $\delta$-function
of the submanifold $\grc_{k/2,n} \subset \gr_{k,2n}$. We will
denote it by $\delta _{\grc}$.
\begin{lemma} Let $k>n$ be even. Then
$\dgrc$ is not orthogonal to the irreducible subspace in
$C(\gr_{k,2n})$ with highest weight $(\underbrace{2,2, \dots,
2}_{2n-k \mbox{ times }}, 0,\dots , 0)$.
\end{lemma}

First let us deduce our statement from this lemma. The cosine
transform commutes with the action of $SO(2n)$ on $C(\gr_{k,2n})$.
Hence by the Schur lemma it acts on each irreducible subspace as a
multiplication by a scalar. Hence an irreducible subspace is
contained in the image of the cosine transform if and only if the
cosine transform on it does not vanish. However by Lemma 2.1.2 the
irreducible subspace with the highest weight vector
$(\underbrace{2,2, \dots, 2}_{2n-k \mbox{ times }}, 0,\dots , 0)$
is contained in the image of $Val^{ev}_{2n,k}$ in $C(\gr_{k,2n})$,
and this image coincides with the image of the cosine transform by
Theorem 1.1.3 of \cite{alesker-bernstein}. Thus it remains to
prove Lemma 2.1.5 to finish the proof of Theorem 2.1.1.

{\bf Proof} of Lemma 2.1.5. First observe that the statement of
the lemma is purely representation theoretical. So replacing each
subspace by its orthogonal complement we may and will assume that
$k\leq n$ (oppositely to our previous assumption on $k$). Under
this assumption it is easier to write down explicit formulas. It
is sufficient to prove that $\dgrc$ is not orthogonal to the
highest weight vector in the relevant irreducible subspace. This
statement will be proven by a computation involving explicit form
of the highest weight vector. First we will write it down
following \cite{strichartz} (see also \cite{grinberg}).

Let $e_{i,j}$ denote $(2n\times 2n)$- matrix which has zeros at
all but one place $(i,j)$ where it has 1. Let us fix a Cartan
subalgebra of $so(2n)$ spanned by $\{C_i\}_{i=1}^n$ where
$$C_i=e_{2i-1,2i}- e_{2i,2i-1},\, i=1, \dots ,n.$$
For any subspace $E\in \gr_{k,2n}$ let us choose an orthonormal
basis $X^1, \dots , X^k$, and let us write its coordinates in the
standard basis in columns of $2n \times k$- matrix:
$$ \left | \begin{array}{ccc}
X_1^1& \dots & X_1^k\\
\multicolumn{3}{c}{\dotfill}\\
 X_{2n}^1& \dots &X_{2n}^k
\end{array}
\right| .$$ Let $X_j$ denote the $j$-th row of this matrix. Let
$A(l)$ be $l\times k$- matrix whose $j$-th row is $X_{2j-1} +
\sqrt{-1} X_{2j}, \, j=1, \dots ,l$. The next lemma was proved in
\cite{strichartz}, Theorem 5 (see also \cite{grinberg}).
\begin{lemma}
Let $k\leq n$. Let $(2m_1, 2m_2, \dots, 2m_k, 0, \dots, 0)$ be the
highest weight of $SO(2n)$ with $m_1\geq m_2 \geq \dots \geq m_k
\geq 0$. The irreducible subspace of $C(\gr _{k,2n})$ with this
highest weight has the highest weight vector of the form
$$ f_{m_1, \dots, m_k}= det[A(1) \cdot A(1)^t]^{m_1 -m_2} \cdot
det[A(2) \cdot A(2)^t]^{m_2 -m_3} \cdot \dots \cdot det[A(k) \cdot
A(k)^t]^{m_k}.$$
\end{lemma}

Recall that we are interested in the highest weight
$(\underbrace{2,2, \dots, 2}_{k \mbox{ times }}, 0,\dots , 0)$.
Hence the highest weight vector $F\in C(\gr _{k,2n})$ has the
form:
$$F :=det[A(k) \cdot A(k)^t].$$
Let us denote for brevity $m:=k/2$ (recall that $m$ is an
integer). We have to show that
$$\int _{M \in \grc _{m,n}} F(M) dM \ne 0.$$
In fact we will show that the function $F$ is non-negative on
$\grc _{m,n}$ and not identically zero.

Let us choose in our hermitian space $\CC^n$ an orthonormal
hermitian basis $e_1, \dots , e_n$. Then in the realization $\RR
^{2n}$ of this space we will choose the basis
$$e_1, e_2, e_3, \dots, e_k ; \sqrt{-1}e_1, -\sqrt{-1}e_2,
\sqrt{-1}e_3,  \dots , -\sqrt{-1}e_k, \mbox{ other vectors}. \eqno
(*)$$ Fix any $E\in \grc _{m,n}$. Let us choose in $E$ an
orthonormal hermitian basis $\xi _1, \dots ,\xi _m$. Then
$$\xi _t =\sum _{j=1}^n z_t^j e_j
=\sum _{j=1}^n (Re z_t^j \cdot e_j + Im z_t^j \cdot (\sqrt{-1}
e_j)),$$ with $z_t^j \in \CC$. Then the vectors $\xi_1,\dots,
\xi_m, \sqrt{-1}\xi_1,\dots,\sqrt{-1}\xi_m$ form a real basis of
$E$. Let us write the coordinates of these vectors with respect to
the basis $(*)$ in columns of the following matrix:
$$ \left[| \begin{array}{ccc|ccc}
Re z_1^1&\dots&Re z_m^1&-Im z_1^1&\dots&-Im z_m^1\\
Re z_1^2&\dots&Re z_m^2&-Im z_1^2&\dots&-Im z_m^2\\
\multicolumn{3}{c}{\dotfill}&\multicolumn{3}{c}{\dotfill}\\
Re z_1^{k-1}&\dots&Re z_m^{k-1}&-Im z_1^{k-1}&\dots&-Im z_m^{k-1}\\
Re z_1^{k}&\dots&Re z_m^{k}  &-Im z_1^{k}  &\dots&-Im z_m^{k}\\
\hline
 Im z_1^1&\dots&Im z_m^1&Re z_1^1&\dots&Re z_m^1\\
-Im z_1^2&\dots&-Im z_m^2&-Re z_1^2&\dots&-Re z_m^2\\
\multicolumn{3}{c}{\dotfill}&\multicolumn{3}{c}{\dotfill}\\
Im z_1^{k-1}&\dots&Im z_m^{k-1}&Re z_1^{k-1}&\dots&Re z_m^{k-1}\\
-Im z_1^k&\dots&-Im z_m^k&-Re z_1^k&\dots&-Re z_m^k\\
\hline
\multicolumn{3}{c}{\dotfill}&\multicolumn{3}{c}{\dotfill}\\
\end{array}
\right].
$$

Now let us write down the $(k\times k)$- matrix $A(k)$. Recall
that the $j$-th row of it is obtained by adding to $(2j-1)$-th row
of the above matrix $i=\sqrt{-1}$ times the $(2j)$-th row of it.
Then we obtain that $A(k)$ is equal to
$$\left[ \begin{array}{ccc|ccc}
Re z_1^1 +i Re z_1^2&\dots&Re z_m^1 +i Re z_m^2&
-Im z_1^1 -i Im z_1^2 &\dots&-Im z_m^1 -i Im z_m^2\\
\multicolumn{3}{c}{\dotfill}&\multicolumn{3}{c}{\dotfill}\\
Re z_1^{k-1} +i Re z_1^k&\dots&Re z_m^{k-1} +i Re z_m^k& -Im
z_1^{k-1} -i Im z_1^k &\dots&-Im z_m^{k-1} -i Im z_m^k\\
 \hline
 Im z_1^1 -iIm z_1^2&\dots&Im z_m^1 -iIm z_m^2& Re z_1^1-iRe
 z_1^2&\dots&Re z_m^1-iRe z_m^2\\
\multicolumn{3}{c}{\dotfill}&\multicolumn{3}{c}{\dotfill}\\
Im z_1^{k-1} -iIm z_1^k&\dots&Im z_m^{k-1} -iIm z_m^k& Re
z_1^{k-1}-iRe z_1^k&\dots&Re z_m^{k-1}-iRe z_m^k
\end{array}
\right].
$$
Let us denote by $A$ the $(m\times m)$-sub-matrix of the above
matrix which stays in the upper left part of it, and by $B$ the
$m\times m$- sub-matrix which stays in the lower left part of it.
Then it is easy to see that
$$ A(k)= \left[
\begin{array}{cc}
A&-\bar B\\
B&\bar A
\end{array}
\right].
$$

Then the function $F$ is equal
$$det[A(k)\cdot A(k)^t]= det[A(k)]^2.$$
Let us show that $det[A(k)] \in \RR$. Indeed
$$\overline {det[A(k)]}= det \left[
\begin{array}{cc}
\bar A&-B\\
\bar B&A
\end{array}
\right]=$$
$$ det \left( \left[
\begin{array}{cc}
0&-1\\
1&0
\end{array}
\right]
 \left[
\begin{array}{cc}
A&-\bar B\\
B&\bar A
\end{array}
\right]
\left[
\begin{array}{cc}
0&-1\\
1&0
\end{array}
\right]\right)= det[A(k)].$$ It remains to show that $F\not \equiv
0$. Take $E_0:=span_{\CC}\{e_1,e_3,e_5,\dots,e_{k-1}\}$. Then
$F(E_0)=1$. \qed

Now we will present another basis in the space of unitarily
invariant valuations. As it was mentioned above this basis is more
convenient to obtain integral geometric formulas for non-convex
sets (see Section 3). Let $\agr _{k,2n}$ denote the Grassmannian
of {\itshape affine} real $k$-dimensional subspaces in $\CC ^n
\simeq \RR ^{2n}$. Let $\agrc _{k,n}$ denote the Grassmannian of
{\itshape affine} complex $k$- dimensional subspaces in $\CC ^n$.
Note that $\agr _{k,2n}$ and $\agrc _{k,n}$ have natural Haar
measures which are unique up to a constant. For every non-negative
integers $p$ and $k$ such that $2p\leq k\leq 2n$  let us introduce
the following valuations:
$$U_{k,p}(K)=\int_{E\in \agrc _{n-p,n}} V_{k-2p}(K\cap E) \cdot
dE.$$ Clearly $U_{k,p}\in Val_{k}^{U(n)}(\CC^n)$.
\begin{proposition}
For any non-negative integers $k,\,p$ satisfying $2p\leq k\leq n$
one has
$$U_{k,p}= c_{n,k,p} \cdot \DD(C_{2n-k,n-p}),$$
where $ c_{n,k,p}$ is a non-zero normalizing constant depending on
$n,\,k$ and $p$ only. Hence the valuations $U_{k,p}$ with $0\leq p
\leq \frac{min\{k, 2n-k\}}{2}$ form a basis of the space
$Val_{k}^{U(n)}(\CC^n)$.
\end{proposition}

{\bf Proof.} Clearly the second statement immediately follows from
the first one and Theorem 2.1.1. First we can rewrite the
definition of $U_{k,p}$ as follows:
$$U_{k,p}(K)= \int_{F\in \grc_{p,n}}dF \cdot \int_{x\in F}dx \cdot
V_{k-2p}(K\cap (x+F^\perp)),$$ where $F^\perp$ denotes the
orthogonal complement of $F$.
 Let us compute the image of $U_{k,p}$ in the space
$C(\gr_{k,n})$ under the imbedding described in Section 1. Fix any
$L\in \gr_{k,n}$. Let $D_L$ denote the unit Euclidean ball inside
$L$. Then by a straightforward elementary computation one can
easily see that for $K=D_L$ the inner integral in the last formula
is equal to $c\cdot |cos(L,F)|$, where $c$ is a normalizing
constant. Hence $$U_{k,p}(D_L)= c\cdot \int _{F\in \grc_{p,n}}dF
\cdot |cos(L,F)|= c\cdot \int _{F\in \grc_{p,n}}dF \cdot
|cos(L^\perp,F^\perp)|=$$
$$c\cdot \int _{E\in \grc_{n-p,n}}dE \cdot
|cos(L^\perp,E)|.\eqno(1)$$ It is easy to see that for any
$M\in\gr_{k,2n}$, and for $2k\leq l$
$$C_{k,l}(D_M)=c'\cdot \int_{E\in \grc_{l,n}}dE \cdot
|cos(M,E)|,\eqno(2)$$ where $c'$ is a normalizing constant.
Clearly (1) and (2) imply the theorem. \qed

\section{Integral geometry in $\CC ^n$.}
Using the classification of unitarily invariant valuations
obtained in the previous section, we will establish new integral
geometric formulas for real submanifolds in $\CC ^n$. Note that
these formulas will be valid not only for convex domains, but for
arbitrary piecewise smooth submanifolds of $\CC ^n$ with corners.

The method to obtain the result for non-convex sets using the
convex case is as follows. First one should guess the correct
formula for the general case. Next one can approximate nicely
piecewise smooth set by polyhedral sets. The last set can be
presented as a finite union of convex polytopes. For each convex
polytope and for each finite intersection of them we can apply the
formulas for the convex case. The final result follows by the
inclusion- exclusion principle. In Subsection 3.1 we obtain new
integral geometric formulas in $\CC^n$. In Subsection 3.2 we
compute explicitly the constants in one of these formulas in the
particular case $n=2$. In Subsection 3.3 we discuss another
example of unitarily invariant valuation, Kazarnovskii's
pseudovolume.
\subsection{General results.} Let us denote by $IU(n)$ the group
of all isometries of $\CC^n$ preserving the complex structure.
(Clearly this group is isomorphic to the semidirect product $\CC^n
\rtimes U(n)$.)

Note also that the intrinsic volumes $V_i$ in a Euclidean space
$\RR^N$ can be defined not only for convex compact domains but
also for compact domains with piecewise smooth boundary (even more
generally, for compact piecewise smooth submanifolds with
corners). For instance for a domain $\Omega$ with smooth boundary
they can be defined as follows: $V_i(\Omega):=\frac{1}{N}
M_{N-1-i}(\partial \Omega)$, where for any hypersurface $\Sigma$
$$M_r(\Sigma):={N-1 \choose r}^{-1} \int_{\Sigma} \{k_{i_1},
\dots,k_{i_r}\} d\sigma,$$ where $\{k_{i_1}, \dots,k_{i_r}\}$
denotes the $r$-th elementary symmetric polynomial in the
principal curvatures $k_{i_1}, \dots,k_{i_r}$, and $d\sigma$ is
the measure induced by the Riemannian metric.

Then we can define the expressions $U_{k,p}(\Omega)$ for $0\leq
2p\leq k\leq 2n$ (for convex compact sets they were defined in
Section 2). The correct generalization is as follows:
$$U_{k,p}(\Omega)=\int_{E\in \agrc_{n-p,n}} V_{k-2p}(\Omega\cap
E)\cdot dE,$$ where we use the above definition of
$V_{k-2p}(\Omega)$.

{\bf Remark.} In fact the expressions $U_{k,p}$ can be defined
also for compact piecewise smooth submanifolds of $\CC^n$ with
corners.

\begin{theorem}
Let $\Omega_1, \, \Omega_2$ be compact domains in $\CC^n$ with
piecewise smooth boundaries such that for every $U\in IU(n)$ the
intersection $\Omega_1\cap U(\Omega _2)$ has finitely many
components. Then
$$\int_{U\in IU(n)} \chi (\Omega_1\cap U(\Omega_2)) dU=
\sum_{k_1+k_2=2n}\sum_{p_1, p_2} \kappa(k_1,k_2,p_1,p_2)
U_{k_1,p_1}(\Omega_1) U_{k_2,p_2}(\Omega_2),$$ where the inner sum
runs over $0\leq p_i\leq k_i/2, \, i=1,2$, and
$\kappa(k_1,k_2,p_1,p_2)$ are certain constants depending on
$n,k_1,k_2,p_1,p_2$ only.
\end{theorem}
\begin{theorem}
Let $\Omega$ be a compact domain in $\CC^n$ with piecewise smooth
boundary. Let $0<q<n, \, 0<2p<k<2q$. Then
$$\int_{E\in \agrc_{q,n}} U_{k,p}(\Omega \cap E)=\sum
_{p=0}^{[k/2]+n-q}\gamma_p \cdot U_{k+2(n-q),p}(\Omega),$$ where
the constants $\gamma_p$ depend only on $n,\,q$, and $p$.
\end{theorem}
Let us denote by $\algr_n$ the (non-compact) Grassmannian of
{\itshape affine} Lagrangian subspaces of $\CC^n$. Clearly it is a
homogeneous space of the group $IU(n)$.

\begin{theorem}
Let $\Omega$ be a compact domain in  $\CC^n$ with piecewise smooth
boundary. Then
$$\int_{\algr (\CC^n)} \chi(E\cap\Omega)dE=
\sum_{p=0}^{[n/2]}\beta_p\cdot U_{n,p}(\Omega),$$ where $\beta_p$
are certain constants depending on $n$ and $p$ only.
\end{theorem}

{\bf Remarks.}1) Theorems 3.1.1 and 3.1.2 are analogs of general
kinematic formulas of Poincar\'e, Chern \cite{chern1},
\cite{chern2} and Federer \cite{federer} (see also \cite{santalo},
especially Ch. 15).

2) These results can be formulated and proved not only for domains
but also for piecewise smooth compact submanifolds in $\CC^n$ with
corners. To do it, consider an $\eps$-neighborhood of this
submanifold for small $\eps >0$. Then apply the above formulas to
this domain. Both sides depend polynomially on $\eps$. Comparing
the lowest degree terms we get the mentioned generalizations. We
do not reproduce here the explicit computations.

3) It would be of interest to compute the constants
$\kappa(k_1,k_2,p_1,p_2)$ and $\beta_p$ in Theorems 3.1.1, 3.1.2,
and 3.1.3 explicitly. We could not do it in general. But we
compute them in the first non-trivial case $n=2$ for Theorem 3.1.3
in the next subsection.

\subsection{Integral geometry in $\CC^2$.}
In this subsection we will compute explicitly the constants in one
of the integral geometric formulas discussed in the previous
subsection in the particular case of $\CC^2$. In order to do these
computations we first recall the classical presentations for the
orthogonal group $SO(4)$ (more precisely for its universal
covering $Spin(4)$) and for the Grassmannian of oriented 2-planes
in $\RR ^4$ which we will denote by $\gro$.

Let us denote the standard complex structure on $\CC^2$ by $i$.
Let us identify $\CC^2$ with the quaternionic space $\HH$ with the
usual anti-commuting complex structures $i,\, j,$ and $k=ij$. Then
clearly $\HH= \CC \oplus j \CC $. Recall that the group of
quaternions with norm equal to 1 acts by left multiplication on
$\HH \equiv \CC^2$ and thus is identified with the group $SU(2)$.

Also we have the isomorphism
$$\Phi: SU(2)\times SU(2)/\{\pm Id\}\tilde{\str} SO(4)$$
defined by
$$\Phi((q_1,q_2))(x)=q_1 x q_2^{-1},$$
where $q_1, \, q_2$ are norm one quaternions. Hence we can and
will identify the group $Spin(4)$ with $SU(2)\times SU(2)$. Let
$E_0 \in \gro$ be the $span_{\RR}\{1,i\}$ with standard
orientation coming from the complex structure. Clearly the
stabilizer of $E_0$ in $Spin(4)=SU(2)\times SU(2)$ is equal to
$T\times T$ where
$$T=\{z\in\CC| |z|=1\}=U(1) \subset SU(2).$$
Hence we have the following presentation of the Grassmannian of
real oriented 2-planes in $\RR^4$ :
$$ \gro =SU(2)/T \times SU(2)/T.$$
However $SU(2)/T\simeq \CC P^1,$ where $ \CC P^1$ is (as usual)
the complex projective line. For  our computations it will be
convenient to identify $\CC P^1$ with the 2-dimensional sphere of
radius $1/2$. Moreover it will be convenient to consider this
sphere $S^2$ in the standard coordinate space $\RR ^3$ with the
center $(1/2,0,0)$. Moreover $E_0\in \gro =S^2\times S^2$ will
correspond to the point $((1,0,0),(1,0,0))$. The following lemma
can be proved by a straightforward computation.
\begin{lemma}
Let $E=(t_1,t_2)\in S^2\times S^2=\gro$. Let $t_i=(x_i,y_i,z_i),\, i=1,2$.
Then $|cos(E,E_0)|=|x_1 +x_2 -1|$.
\end{lemma}
For $\CC^2$ Theorem 2.1.1 says
\begin{proposition} For $0\leq k\leq 4, k\ne 2$ the space
$Val^{U(2)}_{k}(\CC^2)$ is spanned by $V_k$;
\newline
$Val^{U(2)}_{2}(\CC^2)$ is spanned by $V_2$ and by $\phi$, where
$\phi(K)=\int_{\xi\in \textbf{C} P^1} vol_2(Pr_\xi K)d\xi$.
\end{proposition}
Recall that the total measure of $\CC P^1$ is chosen to be equal
to one. Now let us describe the image of the valuation $\phi$ in
$C(\gro)$. Let us denote this image by $f$.
\begin{lemma}
For every $E=(t_1,t_2)\in S^2\times S^2=\gro$ with $t_i=(x_i,y_i,z_i),\, i=1,2$
$$f(E)=vol_2 D_2 \cdot \left((x_2-\frac{1}{2})^2+ \frac{1}{4}\right),$$
where $D_2$ denotes the unit 2-dimensional Euclidean disk.
\end{lemma}
This lemma follows immediately from Lemma 3.2.1 and the fact that
the set of complex lines in $\CC ^2$ is $SU(2)$-orbit of $E_0$.

Let us denote by $LGr_n$ the Grassmannian of Lagrangian subspaces
in $\CC ^n$. Let us define a valuation $\psi \in
Val^{U(2)}_{2}(\CC^2)$ as follows
$$\psi(K)= \int_{F\in LGr_2} vol_2 (Pr_F(K)) dF,$$
where $dF$ is the Haar measure on $LGr_2$ normalized by 1. The
main result of this subsection is
\begin{theorem}
$$\phi +2\psi = V_2.$$
\end{theorem}
{\bf Proof.} Let
$$\tphi :=\frac{1}{vol_2 D_2} \phi =\frac{1}{\pi} \phi,\,
\tpsi :=\frac{1}{vol_2 D_2} \psi =\frac{1}{\pi} \psi.$$
 Let us denote by $\tg$ the image of $\tpsi$ in $C(\gro)$, and
by $\tf$ the image of $\tphi$ in  $C(\gro)$. By Lemma 3.2.3
$$\tf (E) =(x_2- \frac{1}{2} )^2 +\frac{1}{4} \eqno{(**)}$$
for every $E=(t_1,t_2)\in S^2\times S^2= \gro$ with $t_i
=(x_i,y_i,z_i) ,\, i=1,2$. Now let us describe $\tg$.  Thus $\tf$
can and will be considered as a function on the second copy of
$S^2$. Let $E_1=span\{1,j\} \in LGr_2$. It is easy to see that
$E_1= U_0(E_0)$, where $U_0\in SO(4)$ is defined by
$$U_0(x)= \frac{1+k}{\sqrt 2}\cdot x \cdot \overline { \frac{1+k}{\sqrt 2}}$$
for every $x\in \HH =\RR^4$.
Then $$\tg (E)= \int_{F\in LGr_2} |cos(F,E)| dF= \int_{U\in U(2)}
|cos(U(E_1),E)| dU =$$
$$\int_{U\in U(2)} |cos (UU_0 (E_0),E)|dU= \int_{U\in U(2)} |cos(E_0,U_0^{-1} U(E))|dU.$$
However $U(2)= (SU(2)\times U(1))/\{\pm 1\}$, where $(q,\lambda)\in SU(2)\times U(1)$
acts on $x\in \HH$ by $x \mapsto q\cdot x \cdot \lam ^{-1}$. In the formulas below we
will write the action of $\lam \in U(1)$ on $E\in \gro$ from the right: $\lam(E)= E\cdot
\lam ^{-1}$. In this notation the last integral can be rewritten as
$$ \int_{V\in SU(2)}dV \int _{\lam \in U(1)} d \lam \cdot
 |cos(E_0, \frac{1-k}{\sqrt 2} V \cdot E \cdot \lam ^{-1} \cdot
\frac{1+k}{\sqrt 2})|=$$
$$  \int_{V\in SU(2)}dV \int _{\lam \in U(1)} d \lam \cdot
 |cos(E_0, V \cdot E \cdot \lam ^{-1}
\frac{1+k}{\sqrt 2})| =\int_{\lam \in U(1)} d\lam \cdot \tf (E\cdot \lam ^{-1}
 \cdot \frac{1+k}{\sqrt 2}).
$$
By $(**)$ we can write $\tf$ as
$$\tf (E)= h(E) +\frac{1}{3},$$
where $h(E)= (x_2-\frac{1}{2})^2 -\frac{1}{12} $.
The function $h$ on the sphere $S^2$ (of radius $1/2$) has the property
$$\int_{S^2} h=0.$$
We have
$$\tg(E)= \int_{\lam \in U(1)} d\lam \cdot h(E\cdot \lam ^{-1}\cdot \frac{1+k}{\sqrt 2})
+\frac{1}{3}.$$
It is easy to see that
$$h_1(E):= h(E\cdot \frac{1+k}{\sqrt 2})= y_2 ^2 -\frac{1}{12},$$
$-1/2\leq y_2 \leq 1/2$. Clearly $h_1$ is a polynomial of second
degree on sphere $S^2$ such that $\int _{S^2} h_1 =0$. Hence
$\int_{\lam \in U(1)} h_1(E\cdot \lam ^{-1})$ also satisfies these
properties, and moreover it is $U(1)$-invariant. But such a
polynomial is unique up to proportionality, hence $\int_{\lam \in
U(1)} h_1(E\cdot \lam ^{-1}) =c \cdot h(E),$ where $c$ is a
constant. Let us compute it. If subspace $E$ is such that
$x_2=1/2$ then $h(E)=-1/12 $. But
$$h_1(E)= \frac{1}{2\pi} \int_0^{2\pi} d\phi (\frac{cos^2 \phi}{4}- \frac{1}{12})
= 1/24.$$ Hence $c=-1/2$. Thus $\tg(E)=-\frac{h(E)}{2}
+\frac{1}{3}=-\frac{1}{2}(\tf (E)-\frac{1}{3})
+\frac{1}{3}=\frac{1}{2} - \frac{\tf(E)}{2}$. Hence $\tphi +2\tpsi
= \kappa \cdot V_2$, where $\kappa$ is a normalizing constant such
that for the unit 2-dimensional Euclidean disk $D_2$, $\kappa
\cdot V_2(D_2) =1$. Thus we get that $\phi
+2\psi=\frac{\pi}{V_2(D_2)}V_2=\frac{\pi}{vol_2(D_2)}V_2=V_2.$
\qed
\subsection{Kazarnovskii's pseudovolume.} In this subsection we
discuss another example of unitarily invariant translation
invariant continuous valuation which has rather different origin,
namely it comes from complex analysis. We discuss so called
Kazarnovskii's pseudovolume.  The main result of this subsection
is a new formula for Kazarnovskii's pseudovolume in integral
geometric terms. The proof of this result is based on the
classification of unitarily invariant valuations (Theorem 2.1.1).

Now let us recall the definition of Kazarnovskii's pseudovolume
following \cite{kazarnovskii1}, \cite{kazarnovskii2}. Let $\CC ^n$
be Hermitian space with the Hermitian scalar product $(\cdot,
\cdot)$. For a convex compact set $K\in {\cal K}(\CC ^n)$ let us
denote its supporting functional $$h_K(x):= \underset{y\in
K}{\sup}(x,y).$$ For a set $K\in {\cal K}(\CC ^n)$ such that its
supporting functional $h_K$ is smooth on $\CC^n -\{0\}$
Kazarnovskii's pseudovolume $P$ is defined as follows:
$$P(K):=\int _{D} (dd^c h_K)^n ,$$
where $D$ denotes the unit Euclidean ball on $\CC^n$, and $d^c=
I^{-1} \circ d \circ I$ for our complex structure $I$.
\begin{proposition}
Kazarnovskii's pseudovolume $P$ extends by continuity in the
Hausdorff metric to all ${\cal K}(\CC ^n)$. Then $P$ is unitarily
invariant translation invariant continuous valuation homogeneous
of degree $n$.
\end{proposition}
{\bf Proof.} The first part of the proposition (the continuity) is
a standard fact from the theory of plurisubharmonic functions
originally due to Chern-Levine-Nirenberg
\cite{chern-levine-nirenberg} (see also \cite{kazarnovskii1},
\cite{kazarnovskii2}). The unitary invariance, translation
invariance, and the homogeneity of degree $n$ are obvious. The
only thing which remains to prove is that $P$ is a valuation.

Let $A$ be a convex polytope. It was shown by Kazarnovskii
\cite{kazarnovskii1} that
$$P(A)=\kappa \sum _{F} f(F) \gamma (F)vol_n F,$$
where $\kappa$ is a normalizing constant, the sum runs over all
$n$-dimensional faces $F$ of $A$, $\gamma (F)$ is the measure of
the exterior angle of $A$ at $F$, $vol_n F$ denotes the
($n$-dimensional) volume of the face $F$, and $f(F)$ is defined as
follows. Let $D_F$ denote the unit ball in the {\itshape linear}
subspace parallel to the face $F$. Then $f(F)= vol(D_F +I\cdot
D_F)$. It is easy to see from the above formula that $P$
restricted to the class of convex compact polytopes is a
valuation, namely if $A_1, \, A_2, \, A_1\cup A_2$ are convex
compact polytopes then
$$P(A_1\cup A_2)= P(A_1)+P(A_2)-P(A_1\cap A_2).$$
Then it is easy to see that the continuity of $P$ and the
valuation property on the subclass of polytopes imply that $P$ is
a {\itshape weak} valuation on ${\cal K}(\CC^n)$ (this means that
for any real hyperplane $H$ and any $K\in {\cal K}(\CC^n)$ one has
$P(K)=P(K\cap H^+)+P(K\cap H^-)-P(K\cap H)$ where $H^+$ and $H^-$
denote the half-spaces). However it was shown by Groemer
\cite{groemer} that every continuous weak valuation is valuation
(in the usual sense). \qed

The main result of this subsection is as follows.
\begin{theorem}
$$P=\sum_{n/2\leq l\leq n} \alpha _l C_{n,l},$$
where $\alpha _l\in \RR$ are certain constants depending only on
$n$, and $C_{n,l}$ are valuations defined in the previous section.
\end{theorem}

{\bf Remark.} It would be interesting to compute explicitly the
constant $\alpha_l$. We can do it only for $n=2$ by a direct
computation.

{\bf Proof.} The proof follows immediately from Proposition 3.3.1
and Theorem 2.1.1. \qed


\vskip 0.7cm

\begin{thebibliography}{99}
\bibitem{alesker1}
 Alesker, Semyon; On P. McMullen's conjecture on translation invariant valuations.
Adv. Math. 155 (2000), no. 2, 239--263.

\bibitem{alesker2}
Alesker, Semyon; Description of translation invariant valuations
with the solution of P. McMullen's conjecture. Geom. Funct. Anal.
11 (2001), no. 2, 244--272.

\bibitem{alesker-ecm}
Alesker, Semyon; Classification results on valuations on convex
sets. Proceedings of 3-rd ECM, 2000, Barcelona.

\bibitem{alesker-icm}
Alesker, Semyon; Algebraic structures on valuations, their
properties and applications. Proceedings of ICM 2002, Beijing.
\bibitem{alesker-mult}
Alesker, Semyon; The multiplicative structure on polynomial
continuous valuations. Geom. Funct. Anal., to appear. electronic
version in arXiv: math.MG/0301148
\bibitem{alesker-bernstein}
Alesker, Semyon; Bernstein, Joseph; Range characterization of the
cosine transform on higher Grassmannians. Adv. Math., to appear.
electronic version in arXiv: math.MG/0111031

\bibitem{bernstein}
 Bernstein, D. N.; The number of roots of a system of equations.
 Funkcional. Anal. i Prilo\v zen. 9 (1975), no. 3, 1--4.

\bibitem{borho-brylinski}
Borho, Walter; Brylinski, Jean-Luc; Differential operators on
homogeneous spaces. I. Irreducibility of the associated variety
for annihilators of induced modules. Invent. Math. 69 (1982), no.
3, 437--476.

\bibitem{casselman}
Casselman, William; Canonical extensions of Harish-Chandra modules
to representations of $G$. Canad. J. Math. 41 (1989), no. 3,
385--438.

\bibitem{cheeger-muller-schrader}
Cheeger, J.; M\"uller, W.; Schrader, R.; Kinematic and tube
formulas for piecewise linear spaces. Indiana Univ. Math. J. 35
(1986), no. 4, 737--754.
\bibitem{chern1}
 Chern, Shiing-shen; On the kinematic formula in the Euclidean space of $n$ dimensions.
  Amer. J. Math. 74, (1952). 227--236.
\bibitem{chern3}
Chern, Shiing-shen; Geometry of submanifolds in a complex
projective space. 1958 Symposium internacional de topologia
algebraica International symposium on algebraic topology pp.
87--96 Universidad Nacional Aut\'onoma de M\'exico and UNESCO,
Mexico City

\bibitem{chern2}
 Chern, Shiing-shen; On the kinematic formula in integral geometry.
 J. Math. Mech. 16 1966 101--118.
\bibitem{chern-levine-nirenberg}
Chern, S. S.; Levine, Harold I.; Nirenberg, Louis; Intrinsic norms
on a complex manifold. 1969 Global Analysis (Papers in Honor of K.
Kodaira) pp. 119--139 Univ. Tokyo Press, Tokyo.

\bibitem{federer}
Federer, Herbert; Curvature measures. Trans. Amer. Math. Soc. 93,
1959, 418--491.
\bibitem{fu}
 Fu, Joseph H. G.; Kinematic formulas in integral geometry.
 Indiana Univ. Math. J. 39 (1990), no. 4, 1115--1154.

\bibitem{gelfand-graev-rosu}
Gelfand, I. M.; Graev, M. I.; Ro\c su, R.;
 The problem of integral geometry and intertwining operators for a pair of real Grassmannian
manifolds. J. Operator Theory 12 (1984), no. 2, 359--383.
\bibitem{gelfand-kapranov-zelevinsky}
Gelfand, I. M.; Kapranov, M. M.; Zelevinsky, A. V.; Discriminants,
resultants, and multidimensional determinants. Mathematics: Theory
\& Applications. Birkh\"auser Boston, Inc., Boston, MA, 1994.

\bibitem{goodey-howard-reeder}
 Goodey, Paul; Howard, Ralph; Reeder, Mark;
 Processes of flats induced by higher-dimensional processes. III.
 Geom. Dedicata 61 (1996), no. 3,
257--269.

\bibitem{griffiths}
Griffiths, Phillip A.; Complex differential and integral geometry
and curvature integrals associated to singularities of complex
analytic varieties. Duke Math. J. 45 (1978), no. 3, 427--512.

\bibitem{griffiths-harris}
 Griffiths, Phillip; Harris, Joseph; Principles of algebraic geometry.
 Pure and Applied Mathematics.
 Wiley-Interscience [John Wiley and Sons], New York, 1978.
\bibitem{grinberg}
Grinberg, Eric L.; On images of Radon transforms. Duke Math. J. 52
(1985), no. 4, 939--972.

\bibitem{groemer}
Groemer, H.; On the extension of additive functionals on classes
of convex sets. Pacific J. Math. 75 (1978), no. 2, 397--410.

\bibitem{hadwiger-book}
Hadwiger, H.; Vorlesungen \"uber Inhalt, Oberfl\"ache und
Isoperimetrie. (German) Springer-Verlag,
Berlin-G\"ottingen-Heidelberg 1957.

\bibitem{howard}
 Howard, Ralph; The kinematic formula in Riemannian homogeneous spaces.
 Mem. Amer. Math. Soc. 106 (1993), no. 509

\bibitem{howe-lee}
Howe, Roger; Lee, Soo Teck; Degenerate principal series
representations of $ GL _n(\CC)$ and $ GL _n(\RR)$. J. Funct.
Anal. 166 (1999), no. 2, 244--309.

\bibitem{hug-schneider}
Hug, Daniel; Schneider, Rolf; Kinematic and Crofton formulae of
integral geometry: recent variants and extensions. Preprint, 2002.

\bibitem{kang-tasaki1}
Kang, Hong Jae; Tasaki, Hiroyuki; Integral geometry of real
surfaces in the complex projective plane. Geom. Dedicata 90
(2002), 99--106.

\bibitem{kang-tasaki2}
Kang, Hong Jae; Tasaki, Hiroyuki; Integral geometry of real
surfaces in complex projective spaces. Tsukuba J. Math. 25 (2001),
no. 1, 155--164.

\bibitem{kazarnovskii1}
Kazarnovski\u \i , B. Ya.; On zeros of exponential sums. (Russian)
Dokl. Akad. Nauk SSSR 257 (1981), no. 4, 804--808.

\bibitem{kazarnovskii2}
 Kazarnovski\u\i, B. Ya.; Newton polyhedra and roots of systems of exponential sums.
  Funktsional. Anal. i Prilozhen. 18 (1984), no. 4, 40--49, 96.

\bibitem{klain1}
Klain, Daniel; A short proof of Hadwiger's characterization
theorem. Mathematika 42 (1995), no. 2, 329--339.

\bibitem{klain2}
Klain, Daniel; Even valuations on convex bodies. Trans. Amer.
Math. Soc. 352 (2000), no. 1, 71--93.

\bibitem{klain-rota}
Klain, Daniel; Rota, Gian-Carlo; Introduction to geometric
probability. (English. English summary) Lezioni Lincee. [Lincei
Lectures] Cambridge University Press, Cambridge, 1997.

\bibitem{kouchnirenko}
Kouchnirenko, A. G.; The Newton polygon, and Milnor numbers.
Funkcional. Anal. i Prilo\v zen. 9 (1975), no. 1, 74--75.
\bibitem{le-hong-van}
Le Hong Van; Application of integral geometry to minimal surfaces.
Internat. J. Math. 4 (1993), no. 1, 89--111.

\bibitem{matheron1}
Matheron, G.;
 Un th\'eor\`eme d'unicit\'e pour les hyperplans poissoniens. (French)
J. Appl. Probability 11 (1974), 184--189.


\bibitem{mcmullen-euler}
McMullen, Peter; Valuations and Euler-type relations on certain
classes of convex polytopes. Proc. London Math. Soc. (3) 35
(1977), no. 1, 113--135.
\bibitem{mcmullen-survey}
McMullen, Peter; Valuations and dissections.
 Handbook of convex geometry, Vol. A, B, 933--988, North-Holland,
 Amsterdam, 1993.
\bibitem{mcmullen-schneider}
 McMullen, Peter; Schneider, Rolf;
Valuations on convex bodies. Convexity and its applications,
170--247, Birkh\"auser, Basel, 1983.
\bibitem{park}
Park, Heunggi; Kinematic formulas for the real subspaces of
complex space forms of dimension 2 and 3. Ph.D. thesis, 2002.

\bibitem{santalo-paper}
Santal\'o, Luis A.; Integral geometry in Hermitian spaces. Amer.
J. Math. 74, (1952). 423--434.
\bibitem{santalo}
Santal\'o, Luis A.; Integral geometry and geometric probability.
With a foreword by Mark Kac. Encyclopedia of Mathematics and its
Applications, Vol. 1. Addison-Wesley Publishing Co., Reading,
Mass.-London-Amsterdam, 1976.
\bibitem{schneider-book}
 Schneider, Rolf; Convex bodies: the Brunn-Minkowski theory.
 Encyclopedia of Mathematics and its Applications, 44.
 Cambridge University Press, Cambridge, 1993.
 \bibitem{strichartz}
 Strichartz, Robert S.;
The explicit Fourier decomposition of $L^2(SO(n)/SO(n-m))$. Canad.
J. Math. 27 (1975), 294--310.

\bibitem{takeuchi}
Takeuchi, Masaru; Modern spherical functions. Translated from the
1975 Japanese original by Toshinobu Nagura. Translations of
Mathematical Monographs, 135. American Mathematical Society,
Providence, RI, 1994.

\bibitem{tasaki1}
 Tasaki, Hiroyuki; Generalization of K\"ahler angle and integral
 geometry in complex projective spaces.
 Steps in differential geometry (Debrecen, 2000), 349--361, Inst. Math. Inform., Debrecen, 2001.

\bibitem{tasaki2}
Tasaki, Hiroyuki; Integral geometry in complex projective spaces.
Proceedings of the Sixth International Workshop on Differential
Geometry (Taegu, 2001), 23--34, Kyungpook Natl. Univ., Taegu,
2002.

\bibitem{wallach}
 Wallach, Nolan R.;
  Real reductive groups. I.
  Pure and Applied Mathematics, 132. Academic Press, Inc., Boston, MA, 1988.
\end{thebibliography}
\end{document}